\documentclass[10]{article}
\usepackage{amssymb}
\usepackage[mathscr]{euscript}

\newcommand{\R}[1][]{\ensuremath{{\mathbb{R}^{#1}} }}
\newcommand{\C}[1][]{\ensuremath{{\mathbb{C}^{#1}} }}
\newcommand{\Z}[1][]{\ensuremath{{\mathbb{Z}^{#1}} }}
\renewcommand{\P}[1][]{\ensuremath{{\mathbb{P}^{#1}} }}
\newcommand{\D}[1][]{\ensuremath{{\mathbb{D}^{#1}} }}
\newcommand{\RP}[1][]{\R\P[#1]}
\newcommand{\CP}[1][]{\C\P[#1]}
\newcommand{\goth}[1]{\ensuremath{\mathfrak{#1}}}
\newcommand{\fleche}{\mathop{\longrightarrow}}

\renewcommand{\Im}{\mathop{\rm Im}}
\newcommand{\tr}{\mathop{\rm tr}}
\def\<{\langle} \def\>{\rangle}
\newcommand{\1}{\mathrm{1 \hspace{-0.25em} l}}
\newcommand{\bbox}{\normalsize {}%
        \nolinebreak \hfill $\blacksquare$ \medbreak \par}

\newcommand{\trsp}[1]{ {}^t \! #1 }
\newenvironment{matrice}[1]{ \left( \begin{array}{#1}}{ \end{array} \right) }
\newcommand{\bmat}[1]{\begin{matrice}{#1}}
\newcommand{\emat}{\end{matrice}}
\newcommand{\at}[1]{\big|_{#1}}
\newcommand{\g}{\ensuremath{\goth{g}}}
\newcommand{\gC}{\ensuremath{\goth{g}^{\C}}}
\newcommand{\h}{\ensuremath{\goth{h}}}
\newcommand{\m}{\ensuremath{\goth{m}}}
\newcommand{\p}{\ensuremath{\goth{p}}}
\renewcommand{\u}{\ensuremath{\goth{u}}}
\newcommand{\su}{\ensuremath{\goth{su}}}
\renewcommand{\sl}{\ensuremath{\goth{sl}}}

\renewcommand{\phi}{\varphi}
\newcommand{\Ad}{\mathop{\rm Ad}}
\newcommand{\ad}{\mathop{\rm ad}}
\renewcommand{\.}{\cdot}
\newcommand{\lf}{\hfill\break}
\newcommand{\proj}[2]{\left\lfloor #1 \right\rfloor_{#2}}
\newtheorem{theorem}{Theorem}[section]
\newtheorem{definition}[theorem]{Definition}
\newtheorem{proposition}[theorem]{Proposition}
\newtheorem{lemma}[theorem]{Lemma}
\newtheorem{corollary}[theorem]{Corollary}

\newtheorem{example}[theorem]{Example}

\newcommand{\titre}{Hamiltonian stationary Lagrangian surfaces 
in Hermitian symmetric spaces}

\author{
Fr\'ed\'eric H\'elein\thanks{CMLA, ENS de Cachan, 
61 avenue du Pr\'esident Wilson, 94235 Cachan Cedex, France} , 
Pascal Romon\footnotemark[1]\, \thanks{Universit\'e de Marne-la-Vall\'ee, 
5 bd Descartes, Noisy-le-Grand, 77454 Marne-la-Vall\'ee cedex 2, France} 
}
\title{\titre}
\date{}

\begin{document}
\maketitle

\begin{abstract}
This paper is the third of a series on Hamiltonian stationary Lagrangian 
surfaces. We present here the most general theory, valid for any Hermitian 
symmetric target space. Using well-chosen moving frame formalism, we show 
that the equations are equivalent to an integrable system, 
generalizing the $\C^2$ subcase analyzed in~\cite{HR1}. It shares
many features with the harmonic map equation of surfaces into symmetric 
spaces, allowing us 
to develop a theory close to Dorfmeister, Pedit and Wu's, including for 
instance a Weierstrass-type representation. Notice that this 
article encompasses the article mentioned above, although 
much fewer details will be given on that particular flat case.  
\end{abstract}

    \section*{Introduction}

Hamiltonian stationary Lagrangian submanifolds are Lagrangian submanifolds
of a given symplectic manifold $M$ which are critical points  of the
volume functional with respect to a special class of infinitesimal variations
preserving the Lagrangian constraint, namely the class of compactly
supported Hamiltonian vector fields. This makes sense if $M$ is not only 
a symplectic but is also a Riemannian manifold. This is true in particular 
if $M$ is a K\"ahler manifold. The Euler-Lagrange equation has a particularly 
elegant formulation in terms of the so-called \emph{Lagrangian angle} 
$\beta$, a $\R/2\pi\Z$-valued function defined along any Lagrangian 
submanifold. We may think of $\beta$ as part of the Gauss map. A 
Lagrangian surface is Hamiltonian stationary if and only if $\beta$ 
is harmonic.

This problem has been addressed recently 
by some authors like J. Wolfson \cite{Wo}. It offers a nice generalization 
of the minimal surface problem which maybe shares more similarity with the 
constant mean curvature surfaces problem, since the Lagrangian constraint 
could be thought as a replacement for the volume constraint. Along 
these lines, Y. G. Oh studied as a particular solution the Clifford tori
in $\C^2$ or $\CP^2$ and conjectured that these tori actually 
minimize the area among all tori which are isotopic by Hamiltonian 
deformation, a statement which looks like an isoperimetric inequality 
\cite{O}. Other examples of Hamiltonian stationary Lagrangian surfaces 
in $\C^2$ were found by I.~Castro and F. Urbano \cite{CaUr}. Another 
motivation for looking at this problem lies in its similarity to some models 
in incompressible elasticity \cite{Wo,HR1}. Furthermore, the Hamiltonian 
stationary class includes as a subcase the class of special Lagrangian 
manifolds, which are calibrated manifolds introduced by R. Harvey and 
H. B. Lawson \cite{HaL}. These are simply Lagrangian surfaces in 
a Calabi-Aubin-Yau 4-manifold with \emph{constant} Lagrangian angle
function $\beta$. In \cite{ScWo}, R. Schoen and J. Wolfson propose to 
produce special Lagrangian manifolds by first constructing a Hamiltonian 
stationary Lagrangian submanifold using analytical methods and then proving 
that, under some hypotheses, these submanifolds are in fact special 
Lagrangian. These questions are also strongly motivated by string theory 
and branes theory, where special Lagrangian submanifolds play an 
important role.

In \cite{HR1} we considered Hamiltonian stationary Lagrangian 
surfaces in $\R^4\simeq \C^2$ and showed that this problem is completely 
integrable, like the famous KdV equation or, in differential geometry, 
the equation for harmonic maps from a surface to  a symmetric space. 
Our formalism was similar to the harmonic map theory, except that 
the connection form is not of the form $\alpha_\lambda:=
\lambda^{-1}\alpha_1'+\alpha_0 +\lambda \alpha_1''$ but $\alpha_\lambda:=
\lambda^{-2}\alpha_2'+\lambda^{-1}\alpha_{-1}+\alpha_0 +\lambda \alpha_1 
+ \lambda^2\alpha_2''$, where $\lambda$ is a complex (spectral) 
parameter. The  equations being slightly more linear than 
in the harmonic maps context, we were able to simplify the integrable 
system theory and propose a Weierstrass type representation formula 
much simpler than the one constructed by J. Dorfmeister, F. Pedit and 
H. Y. Wu in \cite{DPW}. In \cite{HR2} we further simplified these 
formulas using quaternions and compared them with a similar formula
due to B. G. Konopelchenko, suggesting a hidden structure of spinors. 
Using these formulas, H. Anciaux recently showed estimate towards 
Oh's conjecture \cite{A}.

In the following paper, we shall show that the integrable system structure 
-- with a family of connections of the type $\alpha_\lambda:=
\lambda^{-2}\alpha_2'+\lambda^{-1}\alpha_{-1}+\alpha_0+\lambda \alpha_1 
+ \lambda^2 \alpha_2''$ -- persists if one replaces the ambient space 
by any two-dimensional Hermitian symmetric space. It should be noted 
that the Hermitian symmetric spaces are automatically K\"ahler-Einstein. 
They are of five types, namely : $\R[4] \simeq \C[2]$, \CP[2] and 
its dual the complex hyperbolic space, $\CP[1] \times \CP[1]$ and its dual.
Moreover they share a crucial algebraic property: the existence
of an order four automorphism $\tau$, squaring to the usual symmetry 
$\sigma$, where $M$ is the quotient of all isometries over the fixed point set
of $\sigma$. The specific properties of $\tau$ will be described below.
As an application we show that the theory in \cite{DPW} can be generalized
(fully for compact spaces, partially with local versions for non compact
spaces as in \cite{H}) and that conformal parametrizations of such
surfaces are constructed using holomorphic data. There is no doubt that
it is also possible to build a theory of finite type solutions (in particular
for tori) using our equations, as was done in details in~\cite{HR1}.
However, unlike in the harmonic map setting, the choice of appropriate 
moving frames plays a key role in our theory. We call them 
\emph{Lagrangian framings}. We show how they help to characterize the 
Lagrangian angle and describe their properties in 
great details in section~\ref{symstruct}.

By the same token, we can extend our description to Hamiltonian 
stationary Lagrangian cones in $\C^3$, since their links 
(i.e. intersection with the sphere $S^5$) are stationary Legendrian 
surfaces which project down to Hamiltonian stationary Lagrangian surfaces 
in $\CP^2$ by the Hopf fibration. Note that special Lagrangian cones 
in $\C^3$ with toric links have been constructed by M. Haskins in \cite{Has}.

A last comment: as for harmonic maps, constant mean curvature
surfaces in $\R^3$ or Willmore surfaces \cite{H}, the emergence 
of such a miraculous theory is linked to the existence of conjugate 
families of solutions (which here are obviously obtained by rotating 
the spectral parameter $\lambda$ in $S^1$). There is however 
a novelty with Hamiltonian stationary Lagrangian surfaces, residing 
in its different connection form $\alpha_\lambda$ involved. 
We may furthermore observe that this connection mixes spinor-like 
quantities ($\alpha_{-1}$, $\alpha_1$) and non spinor-like ones 
($\alpha_{-2}$, $\alpha_0$), so we may ask whether there is some 
supersymmetric interpretation of that.

        \section{Moving frames and Lagrangian lifts}

        \subsection{Lagrangian framings of Lagrangian surfaces in a K\"ahler 
        4-manifold }

Let $M$ be a K\"ahler 4-manifold with almost complex structure $J$, 
and $B$ the principal $U(2)$-bundle of unitary frames on $M$ (with 
the obvious $U(2)$ action); define $B' = B/SU(2)$ as the quotient bundle 
(indeed a principal $U(1)$-bundle). Notice that $B'$ is diffeomorphic to 
(dual of) the canonical bundle on $M$. Let $f : L \to M$ be an immersion 
of a surface into $M$. The surface is Lagrangian if $f^* \omega = 0$ 
where $\omega$ is the K\"ahler form; that amounts to saying that for any
$p\in M$, the tangent plane $T_pL$ at $p$ to $L$ is mapped by $J$ to its
(Riemannian) orthogonal complement or, in other words, that any
Riemannian-orthonormal basis $(e_1,e_2)$ of $T_pL$ is actually
Hermitian-orthonormal in $TM$.

Hence for any $p\in L$ the set of all Riemannian orthonormal bases of the
tangent plane $T_pL$ can be described as being $SO(2)_*(e_1,e_2)$ --
i.e. some orbit of the action of $SO(2)$ -- and this set is a subset of 
the set of all Hermitian-orthogonal bases at $p$, i.e. the orbit
$U(2)_*(e_1,e_2)$. In between lies the orbit $SU(2)_*(e_1,e_2)$ which 
precisely interests us.
\begin{definition}      \label{deflagframing}
    A (local) framing of the immersion $f$ is a (local) section of the 
unitary frame bundle $f^*B$.
    A framing $F$ of $f$ is called \emph{Lagrangian}\footnote{this choice 
ought to be compared to Darboux framings in the Riemannian setting,
for a submanifold $N^n$ of $M^m$: a Darboux framing is an orthonormal framing 
such that the $n$ first vectors span the tangent space of $N$, i.e. they 
are in the same class of the frame bundle modulo $SO(n) \times SO(m-n)$. 
As much as Darboux framings give insight to the Riemannian structure, 
Lagrangian framings yield valuable information on the Lagrangian immersion.
See \cite{G} for a description of the moving frame theory.}
    if for any $p\in L$, $F(p)$ is equivalent mod $SU(2)$
    to an orthonormal basis $(e_1(p),e_2(p))$ of $T_pL$.
\end{definition}
Note that there may not exist global sections (Lagrangian or not) of $f^*B$, 
except in coordinate charts. However it is clear that all local Lagrangian
framings lift a single section of the bundle $B'$ by the projection map
$B\to B'$ and because of this uniqueness we obtain a globally
defined section of $B'$.
\begin{definition}\label{deftangentsection}
The image of Lagrangian framings by the quotient map
$B \to B'$ is called the {\em tangent section} of the
bundle $B'$.
\end{definition} 
 
If $M$ is K\"ahler-Einstein, we may describe Lagrangian framings more precisely. 
We follow here the exposition in \cite{Wo} further refined to symmetric 
spaces. Let $K$ denote as is customary the canonical bundle of $M$. Note that
$K$ is canonically isomorphic to $B'$. On the pull-back bundle $f^*K$ with associate
pull-back metric, the first Chern form satisfies
\[ 
f^*c_1 = f^* \textrm{Ric} = R f^*\omega = 0 .
\]
(Here Ric denotes the Ricci form and $R$ the scalar curvature.)
Hence $f^*K$ is flat and we can construct (local) parallel sections.
Any parallel unit section $s_0$ defines a trivialization of $f^*K \simeq
f^*B'$, in which the tangent section is described by a unit complex number 
$e^{i\beta}$. Another way to express that is by evaluating the 
parallel section $s_0$ of $f^*K$ against any orthonormal framing $e_1,e_2$
of $TL$
\[ 
e^{i\beta} = s_0(e_1 , e_2) .
\]
The real valued function $\beta$ is called the Lagrangian angle\footnote{
we differ from Wolfson \cite{Wo} in that $\beta$ is defined mod $2\pi$ instead 
of mod 2.}. A different choice of $s_0$ will only change $\beta$ by a constant.
A Lagrangian framing is a lift with the same Lagrangian angle
as the tangent section. Finally we define the Maslov form as 
$\Theta=\frac{1}{\pi} d\beta$, and it is automatically closed. 
It is useful to have another expression of the Maslov form, namely
$\Theta=\frac{1}{\pi} \iota_H \omega$ where $H$ is the mean curvature 
vector field. Hence Lagrangian surfaces with constant Lagrangian angle are 
automatically minimal and vice-versa. If furthermore the canonical bundle 
on $M$ is flat (as in \C[2] or more generally as
in any Calabi-Aubin-Yau manifold), there is 
a globally defined Lagrangian angle obtained by pulling back a global 
parallel unit section of $K$ (e.g. $dz^1 \wedge dz^2$ in $\C[2]$).
In that particular setting, minimal Lagrangian surfaces are calibrated, 
and thus minimizing; they are called special Lagrangian (see 
\cite{Wo,HaL}).

        \subsection{Symmetric space structure}  \label{symstruct}

Let us now consider a Hermitian symmetric space of the form $M=G/H$, where $G$ 
is the group of unitary transformations of $M$ and $H$ the isotropy group 
at some point henceforth denoted by $p_0$. The symmetry around $p_0$ is 
induced by a group involution $\sigma$ of $G$ with $(G^\sigma)_0 \subset 
H \subset G^\sigma$. Abusing notations, we will also write $\sigma$ for 
$d\sigma(e)$, its differential at identity, acting on $T_{e}G=\g$, and 
similarly, thinking of elements of $G$ as matrices (which indeed they will be) 
we will identify $dg(p)$ with $g$.

The Lie algebra $\g$ (of $G$) is the direct sum of two subspaces 
$\h \oplus \m$, eigenspaces of $\sigma$ with respective 
eigenvalue $+1$ and $-1$ (so that $\h$ is a Lie subalgebra,  
the Lie algebra of $H$). The subspace $\m$ identifies with $T_{p_0}M$.
As a consequence, $\m$ inherits a Hermitian structure. 
The Reader should consult \cite{BRa} for a more thorough introduction to 
homogeneous spaces from a geometer's point of view. 
Fixing a unitary frame $(\epsilon_1, \epsilon_2)$ of $\m$, any element 
$g \in G$ maps that reference frame at $p_0$ to another unitary frame 
$(e_1,e_2)$ at $p = g \. p_0$ through its differential $g$ ($=dg(e)$). 

Suppose now that we can lift the map $f : L \to M$ to $F : L \to G$. 
If, as will often the case, $L$ is a contractible domain, such lifts 
do exist. Thus a choice of $F$ yields a moving frame\footnote{beware 
that the reciprocal is not always true: one cannot always find a lift 
corresponding to a given moving frame; for instance it is true if $H$ 
is four dimensional (as in $\CP[2]$) and false otherwise (as in 
$\CP[1] \times \CP[1]$).}. 
The choice is a priori wide and we have a gauge group $C^\infty(L,H)$ 
acting on lifts. Since we want the lift $F$ to encode some informations about the
first derivatives of $f:L\to M$, we want to 
restrict our possibilities a little bit, and use only Lagrangian lifts, 
i.e. Lagrangian framings in the sense of definition~\ref{deflagframing}.
Lagrangian lifts do exist as will be obvious from the analysis below.
Consider for the moment any lift $F$ and $\alpha = F^{-1}dF$ the associated 
(left) Maurer-Cartan form (the pull-back of the Maurer-Cartan form on $G$).
Using the above symmetric splitting, we write $\alpha = \alpha_{\h} + 
\alpha_{\m}$; the $\m$-valued 1-form $\alpha_\m$ is the pull-back of the 
Maurer-Cartan form of the symmetric space $M$ (as defined in \cite{BRa}); 
it allows us to identify tangent vectors to $L \subset M$ with elements of \m. 
However this identification is gauge dependent: if we replace $F$ with
$Fh^{-1}$ where $h \in C^\infty(L,H)$, $\alpha_{\m}$ changes to
$\Ad h (\alpha_\m)$. Using a reference unitary frame 
$(\epsilon_1,\epsilon_2)$, there exist at any $z \in L$ a unique element
$a(z) \in GL(\m)$ such that
\begin{equation}\label{gl}
\alpha_\m\at{z} = a(z) (\epsilon_1 dx + \epsilon_2 dy).
\end{equation}
Without loss of generality, we may and will assume that $f$ is conformal 
(and $L$ endowed with the induced conformal structure making it a Riemann 
surface; $z=x+iy$ is always a local holomorphic coordinate). Since $(x,y)$ 
are conformal coordinates and $L$ is Lagrangian, (\ref{gl}) can be written as
\begin{equation}        \label{conformalreduction}
\alpha_\m\at{z} = e^{\rho(z)} k(z)(\epsilon_1 dx + \epsilon_2 dy)\; , 
\quad k(z) \in U(\m)
\end{equation}
where $e^{\rho(z)}$ is the conformal factor\footnote{
If $\Ad H$ is big enough -- namely transitive on orthonormal couples -- 
then a suitable gauge change reduces $k(z)$ to the identity; such a lift 
is said \emph{fundamental} (as for instance in $\C^2$, see \cite{HR1}).
However that notion is dependent on 
the coordinate $z$, unlike the notion of \emph{Lagrangian lift} as will 
become clear hereafter.}.

To understand the notion of Lagrangian lift we need to delve deeper 
into the structure of $G$ and $\g$. First, as we shall see in later 
sections, there exists an order four automorphism $\tau$ acting on $G$,
squaring to $\tau^2 = \sigma$. Thus we may
split \gC\/ as $\gC_0 \oplus \gC_2 \oplus  \gC_{-1} \oplus \gC_{1}$, 
where $\gC_{a}$ is the $i^{a}$-eigenspace of $\tau$.
Since $\tau^2 = \sigma$, 
$\h^{\C} = \gC_0 \oplus \gC_2$ and $\m^{\C} = \gC_{-1} \oplus \gC_{1}$. 
Let $G_0$ be the subgroup of $G$ fixed by $\tau$, whose Lie algebra is 
of course $\g_0$. Recall that the adjoint action of $H$ on \m\/
identifies with the action of the linear isotropy 
$H^{*} = \{ dh(p_0) , h \in H \}$ on $T_{p_0} M$
(see \cite{GaHuLa}, Chapter 1); in particular, since all elements of $G$
are unitary transformations of $M$, $\Ad H$ 
identifies with a subgroup of $U(\m)$ (and $\ad \h$ with a 
subalgebra of $\u(\m)$). Moreover, as we shall see later on,
there always is an element $Y \in \h$
such that $\exp( \frac{\pi}{2} \ad Y) = \Ad \exp(\frac{\pi}{2} Y)$
is the complex structure $J$ on $\m$ and the automorphism $\tau$ can be chosen 
so that $\g_2 = \R Y$ (and $\Ad G_0 \subset SU(\m)$).

In order to picture the moving frames, remember that the tangent 
bundle $TN$ is canonically diffeomorphic to the subbundle $[\m]$ of 
$N \times \g$ with fiber $\Ad g (\m)$ over the point $g \. p_0$. The 
moving frame induced by $F$ is $(\Ad F (\epsilon_1), \Ad F(\epsilon_2))$; 
and it is clear now that two lifts $F,F'$ of $f$ define the 
same moving frame mod $SU(2)$ if and only if they are gauge 
equivalent modulo the restricted gauge group $C^\infty(L,G_0)$.
In these local coordinates a unitary framing $(e_1,e_2)\simeq 
\left( e^{-\rho}{\partial f\over \partial x},
e^{-\rho}{\partial f\over \partial y}\right)$
can be identified (in $[\m]$) with
\[
\Ad F\left( e^{-\rho}\alpha_{\m}\left( {\partial \over \partial x}\right) ,
e^{-\rho}\alpha_{\m}\left( {\partial \over \partial y}\right) \right) 
\simeq \Ad F(z)\left( k(z) \epsilon_1,k(z) \epsilon_2\right)
\]
using (\ref{conformalreduction}) -- the later expressions being 
(inconspicuously) independent of the choice of the lift $F$.
The tangent section is the class mod $SU(2)$ of $(e_1,e_2)$. The lift $F$ 
is Lagrangian if and only if its induced moving frame lies in the same class 
as $(e_1,e_2)$. That requires exactly that $k$ be special unitary 
($\det_{\C} k(z) = 1$). That is a very simple condition to check. Such lifts 
do exist on contractible domains since a gauge change $F \to F \exp(-\theta Y)$ 
multiplies $\det_{\C} k$ by $\det_{\C} \Ad \exp(\theta Y) 
= e^{2i \theta}$; we may thus choose $\theta$ to force $k \in SU(\m)$. 
Notice that
\begin{enumerate}
\item[(i)] no integrability condition is involved in this process,
\item[(ii)] $\theta$ is defined mod $\pi$, so that the process can be 
applied to non simply-connected domains, by taking a suitable multiple 
covering (namely double covering for each generator).
\end{enumerate}

How does one read the Lagrangian angle now ? First we have to build 
a parallel lift $\bar{F}$; recall that the covariant derivative is simply 
the flat differentiation followed by the projection on the bundle $[\m]$, 
that we denote here by $\proj{\dots }{[\m]}$, then we have for $i=1,2$
(denoting $\bar{\alpha} = \bar{F}^{-1}d\bar{F}$)
\begin{eqnarray*} 
\nabla \bar{e}_i &=& \proj{d (\Ad \bar{F} (\epsilon_i))}{[\m]}
= \proj{\Ad \bar{F} ([\bar{\alpha},\epsilon_i])}{[\m]}
= \Ad \bar{F} ( \proj{([\bar{\alpha},\epsilon_i])}{\m} )
\\
&=& \Ad \bar{F} ([\bar{\alpha}_{\h},\epsilon_i])
= [ \Ad \bar{F}(\bar{\alpha}_{\h}), \bar{e}_i ]
\end{eqnarray*}
so that the unitary frame $\bar{s}:=(\overline{e}_1,\overline{e}_2)$
varies according to 
$\nabla \bar{s} = \zeta(s)$ with
$\zeta: V\mapsto [\Ad \bar{F} (\bar{\alpha}_{\h}),V]$. 
The class modulo $SU(2)$ of $s$ varies following the class modulo $\su(2)$ 
of $\zeta$, i.e. it is the projection of $\zeta$ on the diagonal 
(central) part $\goth{z}(\Ad \bar{F}(\m))$ in the Lie algebra
decomposition $\u(\Ad \bar{F}(\m)) = 
\goth{z}(\Ad \bar{F}(\m)) \oplus \su(\Ad \bar{F}(\m))$. Since
\[
\proj{\zeta}{\goth{z}(\Ad \bar{F}(\m))} = \Ad \bar{F} 
\proj{\Ad \bar{F}^{-1} \zeta \Ad \bar{F}}{\goth{z}(\m)} \Ad \bar{F}^{-1}
= \Ad \bar{F} \proj{\ad \bar{\alpha}_{\h}}{\goth{z}(\m)} \Ad \bar{F}^{-1}, 
\]
$\bar{F}$ is flat if and only if $\proj{\ad \bar{\alpha}_{\h}}{\goth{z}(\m)} = 0$, 
that is $\bar{\alpha}_{2} = 0$, because the adjoint representation maps 
$\g_2$ to $\goth{z}(\m)$ and $\g_0$ into $\su(\m)$. We now have four ways 
of reading the Lagrangian angle:
\begin{enumerate} 
\item $e^{i\beta}$ is the complex determinant of the tangent frame
$(\Ad \bar{F}(z) ( \bar{k}(z) \epsilon_1)$,
$\Ad \bar{F}(z) ( \bar{k}(z) \epsilon_2))$ in the parallel frame 
$(\Ad \bar{F}(z) (\epsilon_1)$,$\Ad \bar{F}(z) ( \epsilon_2))$, namely 
is the complex determinant of $\bar{k}$ for any parallel framing 
of $f^*B'$,
\item given a Lagrangian lift $F$, and $h = \bar{F}^{-1}F$ the
($H$-valued) gauge change, then $e^{i\beta} = \det_{\C} \Ad h$,
\item in the previous gauge change, write $h$ as a (commutative) product 
$h_2 h_0$ in $G_0 G_2$, which is unique up to sign ($G_2 = \exp \g_2$); 
then $h_2 = \exp(\beta Y/2)$,
\item finally, for a Lagrangian lift, it follows from the previous
characterization (3) that
\begin{eqnarray*}
\frac{d\beta }{2} Y &=& h_2^{-1}dh_2 
= \proj{\Ad h_0(h^{-1}dh)-dh_0 h_0^{-1}}{\g_2}
\\
&=& \proj{h^{-1}dh}{\g_2}
= \proj{\alpha - \Ad F^{-1}(\bar{\alpha})}{\g_2} = \alpha_2.
\end{eqnarray*}
\end{enumerate}
In the following we shall actually exploit the last characterization.
We conclude by noting that the Maslov form is exactly $\frac{2}{\pi}$ 
times the $Y$ component of $\alpha_2$. \\

So far we have not used at all the underlying complex structure of $L$. 
But decomposing the Maurer-Cartan form $\alpha$ into its $(1,0)$ and $(0,1)$ 
parts $\alpha'$ and $\alpha''$ respectively, 
yields an interesting condition for a lift $F$ to be Lagrangian:
\begin{proposition}  \label{Lagrangianispprim}
A lift is Lagrangian if and only if the Maurer-Cartan form
$\alpha = \alpha_2 + \alpha_0 + \alpha_{-1} + \alpha_{1}$
satisfies $\alpha''_{-1} = 0$ (which by reality assumption implies 
$\alpha'_1=0$).
\end{proposition}
\emph{Proof.}
The equivalence rests upon the following simple fact: $\gC_{-1}$
is \emph{exactly} the orbit under $\R_+ \times SU(\m)$ of the vector 
$\epsilon = \frac{1}{2}(\epsilon_1 -i\epsilon_2) \in \m^{\C}$
(for any choice of Hermitian basis $(\epsilon_1,\epsilon_2)$).
A proof is given in~\cite{HR1}. Now writing $\alpha_{\m} = e^\rho
k (\epsilon dz + \bar{\epsilon}d\bar{z})$, we see that $\alpha'_{\m} = 
e^\rho k \epsilon dz$ belongs to $\gC_{-1}$ if and only if $k$ is in 
$SU(\m)$. \bbox
This type of condition should be compared with \emph{primitivity} 
conditions \cite{BP} or $\omega$-maps \cite{Hi}. However our 
``partial primitivity'' differs in that it is a first order requirement
on the immersion (namely that of being Lagrangian) and not an Euler-Lagrange 
condition, while primitive maps and $\omega$-maps are automatically harmonic.
Finally we remark that the assumption of conformality is not crucial 
until the very final step; indeed the concept of Lagrangian lift makes sense 
in any dimension.

        \subsection{Lagrangian surfaces in $\CP^2$ and Lagrangian cones 
        in $\C^3$}    \label{cones}

We concern ourselves here with Lagrangian three dimensional cones 
centered at the origin in \C[3], and in particular with their intersection
with the unit sphere $S^5$, which is called the \emph{link}. We assume 
henceforward that the cone is regular, i.e. the link $M$ is a 
connected submanifold. Notice that our analysis applies as well to regular 
conical singularities in complex three dimensional manifolds.
We recall here the known correspondence between 
Lagrangian cones and Lagrangian surfaces in \CP[2] (see for instance~\cite{Re1,Re2}).

Lagrangian cones in $\C^3$ may be locally described by Lagrangian 
surfaces in $\CP^2$. This is done by the canonical projection map
$P:\C^3\setminus \{0\}\to \CP^2$. More precisely to any Lagrangian 
cone $\Sigma^3$ in $\C^3$ we may associate the Lagrangian surface 
$L:=P(\Sigma^3)$ in $\C^3$. Conversely this surface $L$ describes 
completely $\Sigma^3$ up to a rigid motion in the sense that the set of
Lagrangian cones of $\C^3$ which are mapped to $L$ by
$P$ is $\{ e^{i\alpha}\Sigma^3/\alpha\in \R/2\pi\Z\}$.
This follows from the following picture. Let
$\langle.,.\rangle _H = \langle.,.\rangle _E -i\omega(.,.)$ be the standard
Hermitian product in $\C^3$ ($\langle.,.\rangle _E$ is the Euclidean scalar
product and $\omega(.,.)$ the symplectic form). Let 
$$
S^5:=\{ Z\in \C^3/\langle Z,Z\rangle _H=1\}
$$
be the unit sphere in $\C^3$, and ${\cal L}:=\Sigma\cap S^5$ its link,
which fully determines $\Sigma$.
Now at any point $p\in {\cal L}\subset \Sigma^3$ choose an
orthonormal frame $(e_1,e_2,e_3)$ of $T_p\Sigma^3$ such that $e_3=p$.
Then $(e_1,e_2)$ is an orthonormal basis of $T_p{\cal L}$.
The Lagrangian constraint on $\Sigma$ at the point $p$ can be stated as
$iT_pM\perp T_pM$; this is equivalent to the two conditions
\begin{enumerate}
\item[\bf a)] $T_p{\cal L}\perp ip$ and 
\item[\bf b)] $ie_2\perp e_1$.
\end{enumerate}
The first condition {\bf a)} means
that $T_p{\cal L}$ is contained in $\Pi_p$, the 4-dimensional subspace
of $T_pS^5$ orthogonal to  $ip$. The collection $\Pi:=(\Pi_p)_{p\in S^5}$
forms a distribution on $S^5$, orthogonal to the fibers of the Hopf 
fibration ${\cal H}:S^5\to \CP^2$, and therefore named the {\em horizontal} 
distribution. Notice that each $\Pi_p$ can also be identified with 
the orthogonal subspace to $p$ in $\C^3$ for the Hermitian product. 
The restriction of $\langle.,.\rangle _H$ to $\Pi_p$ is also Hermitian, 
which implies that the restriction of $\omega$ to $\Pi_p$ is symplectic 
(non degenerate). Now the second condition {\bf b)} just means that
$T_p {\cal L}$ is a Lagrangian plane in $\Pi_p$. These two conditions
actually ensure that ${\cal L}$ is a {\em Legendrian} submanifold of the
contact manifold $(S^5,\Pi,\omega_{|\Pi})$.

The horizontality condition {\bf a)} on ${\cal L}$, $T_p{\cal L}\subset \Pi_p$,
$\forall p$ can be translated into the property that the restriction of
the Hopf fibration ${\cal H}$ to ${\cal L}$ is an isometric covering map (and locally
an isometric diffeomorphism onto its image). The second condition {\bf b)} is
then contained in the property that ${\cal H}({\cal L})$ is a Lagrangian
submanifold of $\CP^2$.

Conversely let us start from a contractible Lagrangian $L\subset \CP^2$
and let us try to construct first a Legendrian link lifting ${\cal L}$ in $S^5$ 
and then a Lagrangian cone in $\C^3$. We use the canonical connection $\nabla$
on the Hopf fibration defined for any section
$s:\CP^2\to S^5$, any point $\underline{p}\in \CP^2$ and
any $X\in T_{\underline{p}}\CP^2$ by
$$
\nabla _Xs(\underline{p}) :=
\langle ds_{\underline{p}}(X),s(\underline{p})\rangle_H s(\underline{p}).
$$
The curvature of $\nabla$ is given by
$\nabla_X\nabla_Y s - \nabla_Y\nabla_X s - \nabla_{[X,Y]} s =
2is^{\star}\omega_{\C^3}(X,Y)s$ and hence vanishes along any Lagrangian 
surface. It follows that, if $L$ is Lagrangian and if $f$ denotes
the immersion mapping $L\subset \CP^2$, then $f^{\star}{\cal H}$,
the pull-back of the Hopf bundle by $f$ is flat. Since $L$ is also 
contractible we can construct a flat section $s:L\to S^5$ (unique up 
to the choice of the value of $s$ at one point) and ${\cal L}:=s(L)$ is just 
the Legendrian lift of $L$ that we were looking for. Now this Legendrian 
surface spans a Lagrangian cone $\Sigma^3$ in $\C^3$.\\

Furthermore we observe that $\Sigma^3$ is Hamiltonian stationary if and only if
$L$ is so. We denote $\Theta:=dz^1\wedge dz^2\wedge dz^3$ the 3-form which 
helps us to define the Lagrangian angle function on $\Sigma^3$ by
$e^{i\beta}= \Theta(e_1,e_2,e_3)$ for any orthonormal basis $(e_1,e_2,e_3)$
of $T_p\Sigma^3$. The cone $\Sigma^3$ is Hamiltonian stationary if
and only if $\Delta_{\Sigma^3}\beta=0$. Since $\beta$ obviously does not 
depend on the radius $r:=\langle p,p\rangle _H^{1/2}$, this condition is 
equivalent to $\Delta_{\cal L}\beta=0$ along ${\cal L}$.
Now we need a parallel section $\theta$ of $f^{\star}K$, where $K$ is 
the canonical bundle of $\CP^2$ and $f$ denotes the immersion $L\subset \CP^2$. 
A very simple construction is the following: let $s:L\to {\cal L}$ be 
the lift mapping, it is a parallel section of $f^{\star}{\cal H}$ 
- because ${\cal L}$ is Legendrian - and hence the 2-form 
$\theta:= s^{\star}(\iota_s\Theta)$ is parallel. We can thus use $\theta$
to construct the Lagrangian angle along $L$. Let us denote $\underline{p}$  
a point in $L$ and $p=s(\underline{p})$ its lift in ${\cal L}$. The Lagrangian
angle $\underline{\beta}$ at $\underline{p}$ is defined by
$e^{i\underline{\beta}} = \theta(\underline{e}_1,\underline{e}_2)$,
where $(\underline{e}_1,\underline{e}_2)$ is an orthonormal basis of
$T_{\underline{p}}L$.
But since $(p,s_{\star}\underline{e}_1,s_{\star}\underline{e}_2)$
is an orthonormal basis of $T_p\Sigma^3$,
$\theta(\underline{e}_1,\underline{e}_2) = \Theta(s(\underline{p}),
s_{\star}\underline{e}_1,s_{\star}\underline{e}_2) = e^{i\beta}$
and hence the two Lagrangian angles coincide. It follows that the harmonicity
condition on $\beta$ along ${\cal L}$ is equivalent to the harmonicity condition
on $\underline{\beta}$ along $L$, i.e. the condition that $L$ is Hamiltonian
stationary. 

Another presentation of this relationship between Lagrangian surfaces in 
$\CP[2]$ and Lagrangian cones in $\C[3]$ from the moving frame point 
of view will be found in \S\ref{cones_frames}.

        \section{Loop group formulation}
        
        \subsection{A general formulation of the problem using a family 
        of curvature free connections}
        
In order to yield a well-defined Lagrangian immersion $f$ on the 
contractible domain $L$, the Maurer-Cartan form $\alpha$ of a Lagrangian 
lift $F$ needs only to satisfy the closedness condition (also called 
zero curvature equation, when thinking of $d+\alpha$ as bundle connection): 
$d\alpha + \frac{1}{2}[\alpha \wedge \alpha] = 0$. This equation splits along 
the eigenspaces to give:
\[
\left\{ \begin{array}{l}
d\alpha_2 + [\alpha_0 \wedge \alpha_2] 
+ \frac{1}{2}[\alpha_{-1} \wedge \alpha_{-1}]
+ \frac{1}{2}[\alpha_1 \wedge \alpha_1] = 0
\\
d\alpha_0 + \frac{1}{2}[\alpha_0 \wedge \alpha_0]
+\frac{1}{2}[\alpha_2 \wedge \alpha_2]+ [\alpha_{-1} \wedge \alpha_1] = 0
\\
d\alpha_{-1} + [\alpha_0 \wedge \alpha_{-1}]+[\alpha_2 \wedge \alpha_1]=0
\\
d\alpha_1 + [\alpha_0 \wedge \alpha_1]+[\alpha_2 \wedge \alpha_{-1}]=0
\end{array} \right.
\]     
However due to the commutation relations $[\g_2,\g_2]=[\g_2,\g_0]=0$, 
and to the properties of our Lagrangian lift, we simplify that to:
\begin{equation} \label{flatness}
\left\{ \begin{array}{l}
d\alpha_2 = 0
\\
d\alpha_0 + \frac{1}{2}[\alpha_0 \wedge \alpha_0]
+ [\alpha'_{-1} \wedge \alpha''_1] = 0
\\
d\alpha'_{-1} + [\alpha_0 \wedge \alpha'_{-1}]+[\alpha_2 \wedge \alpha''_1]=0
\\
d\alpha''_1 + [\alpha_0 \wedge \alpha''_1]+[\alpha_2 \wedge \alpha'_{-1}]=0
\end{array} \right.
\end{equation}
As noted in \cite{CM,Wo}, a Lagrangian surface is Hamiltonian stationary, 
i.e. stationary with respect to Hamiltonian deformations, if 
and only if $\beta$ is harmonic, that is $\Theta$ is coclosed, in 
other words $d \star \alpha_2 = 0$. Furthermore this Lagrangian surface 
is minimal\footnote{special Lagrangian if $M = \C^2$.} if and only if
$\alpha_2 = 0$.
A now classical trick allows us to join these two differential equations
into one zero curvature equation, though formulated on a loop algebra.
\begin{proposition}\label{characterization}
On a simply connected domain $L$, the \g-valued one-form $\alpha$ 
is the Maurer-Cartan of a weakly conformal Lagrangian immersion if and 
only if it is flat: $d\alpha + \frac{1}{2}[\alpha \wedge \alpha]=0$
and partially primitive: $\alpha''_{-1}=\alpha'_1=0$.
Furthermore the immersion is Hamiltonian stationary if and only if
the \emph{extended Maurer-Cartan form} 
$\alpha_\lambda = \lambda^{-2} \alpha'_2 + \lambda^{-1}\alpha'_{-1}
+ \alpha_0 + \lambda \alpha''_1 + \lambda^2 \alpha''_2$ is flat
for all $\lambda \in \C^*$ (or $S^1$):
\begin{equation} \label{flatloop}
d\alpha_\lambda + \frac{1}{2}[\alpha_\lambda \wedge \alpha_\lambda]=0
\end{equation}
with the minimal\footnotemark[\value{footnote}] subcase characterized 
by $\alpha_2 = 0$.
\end{proposition}
\emph{Proof.}
Just check that the only additional condition induced 
by~(\ref{flatloop}) is the coclosedness condition on $\alpha_2$. \bbox

Equation (\ref{flatloop}) is a curvature free condition, a.k.a. 
the compatibility condition for the existence for any $\lambda\in S^1$ of 
a solution $F_\lambda:L\to G$ of the equation
\begin{equation}\label{integration}
dF_\lambda = F_\lambda\alpha_\lambda.
\end{equation}
Moreover $F_\lambda$ is unique provided that we know its value at some
point $p_0\in L$. We may choose for instance $F_\lambda(p_0) = \1$.
Is is not difficult to realize that each $F_\lambda$ is a lift of
a Hamiltonian stationary Lagrangian surface like $F$. We deduce 
the following:
\begin{corollary}       \label{family}
Local Hamiltonian stationary Lagrangian surfaces come in families
parametrized by $\lambda \in S^1$. For global surfaces there usually are
period problems.
\end{corollary}

        \subsection{Loop groups}

In the light of Proposition~\ref{characterization} and Corollary~\ref{family} 
above it is natural to introduce the following loop groups 
- following now classical techniques since \cite{SW} (based on ideas
which come back to \cite{SaSa}) and \cite{U}. 
We define the set of maps from $S^1$ to $G$:
\[
\Lambda G:= \{S^1\ni \lambda\mapsto \phi_\lambda\in G\}.
\]
We assume that these maps are bounded in the $H^s$ topology, with
$s>1/2$, where using the Fourier decomposition
$\phi_\lambda=\sum_{k\in \Z}\hat{\phi}_k\lambda^k$,
the $H^s$ norm is defined by
$||\phi_\lambda||_s^2:= \sum_{k\in \Z}|k|^{2s}|\hat{\phi}_k|^2$.
Then $\Lambda G$ is a group for the composition law
$[\lambda\mapsto \phi_\lambda].
[\lambda\mapsto \psi_\lambda] = 
[\lambda\mapsto \phi_\lambda\psi_\lambda]$.
We also consider the twisted loop (sub)group
\[\Lambda G_{\tau}:= \{S^1\ni \lambda\mapsto \phi_\lambda\in G
/ \tau(\phi_\lambda) = \phi_{i\lambda}\}.
\]
These loop groups have Lie algebras which are respectively
\[
\Lambda\g := \{S^1\ni \lambda\mapsto \xi_\lambda\in \g \}
\]
and
\[
\Lambda\g_{\tau}:= \{S^1\ni \lambda\mapsto \xi_\lambda\in \g 
/ \tau(\xi_\lambda) = \xi_{i\lambda}\}.
\]

A key observation is that $\alpha_\lambda$ can be seen as a 1-form with 
values in $\Lambda\g_{\tau}$. More precisely, ``partially primitive'' 
extended 1-forms are exactly the 1-forms $\alpha_\lambda$ with values 
in $\Lambda\g_{\tau}$ such that $\lambda^2 \alpha_\lambda$ has a limit 
when $\lambda$ goes to zero. Similarly, choosing $F_\lambda(p_0) = \1$, 
the family of maps $F_\lambda:L\to G$ 
solution of (\ref{integration}) can rather be viewed as a map into 
$\Lambda G_{\tau}$, called the \emph{extended lift} of $f$.
Such loop groups have already been considered in \cite{BP,DPW} in the context
of harmonic maps into a symmetric space or in \cite{HR1} for Hamiltonian
stationary Lagrangian surfaces in $\C^2$.

        \section{Weierstrass type representations}
        
As in \cite{HR1}, the above loop formulation of the Hamiltonian
stationary Lagrangian surface problem opens the gate to the use of various
constructions of solutions to this problem, using completely integrable systems.
As an illustration, we will present here Weierstrass representations in the
spirit of J. Dorfmeister, F. Pedit and H.Y. Wu \cite{DPW}.

        \subsection{Loop groups decompositions}
        
At the base of this construction is the idea of Iwasawa decomposition.
For instance we shall need such a property for $G_0$, the subgroup of $G$ 
fixed by $\tau$, namely: there exists a solvable Lie subgroup
$B_{G_0}$ of $G_0^{\C}$ such that the following mapping
$$
\begin{array}{ccc}
G_0\times B_{G_0} & \fleche & G_0^{\C}\\
(g,b) & \longmapsto & gb
\end{array}
$$
is a diffeomorphism. We summarize by $G_0^{\C} = G_0B_{G_0}$ this property is
named {\em Iwasawa decomposition}. But we actually need more: an infinite dimensional
extension of this property to loop groups.

For that purpose we need to introduce the complexified versions of the above
loop groups, obtained by replacing $G$ by its complexification $G^{\C}$:
$$
\Lambda G^{\C}:= 
\{S^1\ni \lambda\mapsto \phi_\lambda\in G^{\C}\},
$$
$$
\Lambda G^{\C}_{\tau}:= \{S^1\ni \lambda\mapsto \phi_\lambda
\in G^{\C}/ \tau(\phi_\lambda) = \phi_{i\lambda}\}
$$
and their Lie algebras $\Lambda\g ^{\C}$ and 
$\Lambda\g ^{\C}_{\tau}$. And we also introduce the subgroups
$$
\Lambda^+G^{\C}_{\tau}:= \{[\lambda\mapsto
\phi_\lambda]\in \Lambda G^{\C}_{\tau}
\hbox{ extending holomorphically in the disk }D^2\},
$$
$$
\Lambda^+_{B_{G_0}}G^{\C}_{\tau}:= \{[\lambda\mapsto
\phi_\lambda]\in L^+G^{\C}_{\tau}/
\phi_0 \in B_{G_0}\},
$$
$$
\Lambda^-_{\star}G^{\C}_{\tau}:= \{[\lambda\mapsto
\phi_\lambda]\in \Lambda G^{\C}_{\tau}
\hbox{ extending holomorphically in } S^2 \backslash D^2
\hbox{ and } \phi_{\infty}=1\}.
$$
The two main tools are the following Lemmas, which are proved in \cite{DPW} 
(the proofs are based on Theorems in \cite{PrS}). 
\begin{lemma}\label{iwasawa}
Assume that $G$ is a compact Lie group. Let $\tau:G\to G$ be an order four
automorphism of $G$ and let $G_0$ be the subgroup of $G$ fixed by $\tau$. 
Suppose that the Iwasawa decomposition $G_0^{\C}= G_0B_{G_0}$ holds. 
Then the mapping
$$\begin{array}{ccc}
\Lambda G_{\tau}\times \Lambda^+_{B_{G_0}}G^{\C}_{\tau} & \fleche
& \Lambda G^{\C}_{\tau}\\
(g_\lambda,b_\lambda) & \longmapsto & g_\lambda b_\lambda
\end{array}$$
is a diffeomorphism. We denote by
$\Lambda G^{\C}_{\tau}=\Lambda G_{\tau}.
\Lambda^+_{B_{G_0}}G^{\C}_{\tau}$
this property.
\end{lemma}
\begin{lemma}\label{riemann-hilbert}
Assume that $G$ is a semisimple Lie group. Then there exists a dense open 
subset ${\cal C}$ of the connected component of the identity of 
$\Lambda G^{\C}_{\tau}$, called the {\em big cell}, such that the mapping
\[
\begin{array}{ccc}
\Lambda^-_{\star}G^{\C}_{\tau}\times \Lambda^+G^{\C}_{\tau} 
& \fleche & {\cal C}
\\
(\phi^-_\lambda,\phi^+_\lambda) & \longmapsto 
& \phi^-_\lambda\phi^+_\lambda
\end{array}
\]
is a diffeomorphism. We denote by
${\cal C}=\Lambda^-_{\star}G^{\C}_{\tau}.\Lambda^+G^{\C}_{\tau}$
this property.
\end{lemma}

In some cases in this paper these results do not apply directly, either because 
the isometry group $G$ is not compact (for $\C^2$, $\C\D^2$ 
or the dual of $\CP^1\times \CP^1$) or because this group is not 
semisimple (in the case of $\C^2$). However it is possible to extend 
the above Lemmas to these situations in two ways:
\begin{itemize}
\item for $G = U(2)\ltimes \C^2$, the properties stated in Lemmas 
\ref{iwasawa} and \ref{riemann-hilbert} are true; it was proved in
\cite{HR1} by a direct construction.
\item in all cases, in particular when $G$ is not compact or not semi-simple, 
local versions of Lemmas \ref{iwasawa} and \ref{riemann-hilbert} can 
be proven. In these versions, one just need to replace the loop groups by 
a neighborhood of the identity. The proof of these results uses the inverse 
mapping theorem as in \cite{H}.
\end{itemize}

        \subsection{Solutions in terms of holomorphic data}

Conformal immersions of Hamiltonian stationary Lagrangian surfaces are in 
correspondence with holomorphic data as defined below. We first denote
$$
\Lambda_{-2,\infty}\g ^{\C}_{\tau}:= \{[\lambda\mapsto \xi_\lambda]
\in \Lambda\g ^{\C}_{\tau}/
\xi_\lambda = \sum_{k=-2}^{\infty}\hat{\xi}_k\lambda^k\}.
$$
\begin{definition}
The set of {\em holomorphic potentials}, denoted ${\cal H}_{-2,\infty}(L)$, is 
the set of holomorphic 1-forms on $L$ with values in 
$L_{-2,\infty}\g ^{\C}_{\tau}$. So any
form $\mu_\lambda$ in ${\cal H}_{-2,\infty}(L)$ has the expression
$$
\mu_\lambda = \sum_{k=-2}^{\infty}\hat{\mu}_k\lambda^k
= \sum_{k=-2}^{\infty}\hat{\xi}_k(z) \lambda^k dz,
$$
where $\forall z, \; \sum_{k=-2}^{\infty}\hat{\xi}_k(z)\lambda^k\in
\Lambda_{-2,\infty}\g ^{\C}_{\tau}$.
\end{definition}
\begin{lemma}\label{holopotential}
Let $F_\lambda:L\to \Lambda G_{\tau}$ be the extended lift 
of a (conformal) Hamiltonian stationary Lagrangian immersion and assume that 
$L$ is contractible. Then
\begin{itemize}
\item there exist a holomorphic map $H_\lambda:L\to 
\Lambda G^{\C}_{\tau}$ and a map $B_\lambda:L\to 
\Lambda^+_{B_{G_0}}G^{\C}_{\tau}$ such that 
$F_\lambda = H_\lambda B_\lambda$.
\item the Maurer-Cartan form $\mu_\lambda:= (H_\lambda)^{-1}dH_\lambda$ 
is a holomorphic potential.
\end{itemize}
\end{lemma}
\emph{Proof. } (see \cite{DPW} for details) The existence of $H_\lambda$ and 
$B_\lambda$ relies on solving the equation
\[
0 = {\partial (F_\lambda(B_\lambda)^{-1})\over \partial \overline{z}} =
F_\lambda\left( \alpha_\lambda\left( {\partial \over\partial \overline{z}}
\right)
- (B_\lambda)^{-1}{\partial B_\lambda\over \partial \overline{z}}\right) 
(B_\lambda)^{-1},
\]
which is equivalent to
$$
{\partial B_\lambda\over \partial \overline{z}} = B_\lambda
(\alpha_0 + \lambda\alpha_1 + \lambda^2\alpha_2)\left(
{\partial \over\partial \overline{z}}\right) ,
$$
with the constraint that $B_\lambda$ takes values in $L^+_{B_{G_0}}G^{\C}_{\tau}$.
The existence of a solution is first obtained locally, then we can glue local
solutions into a global one. This proves the first assertion. Now
we write
\[
(H_\lambda)^{-1}dH_\lambda = B_\lambda(\alpha_\lambda
-(B_\lambda)^{-1}dB_\lambda)(B_\lambda)^{-1},$$
and using the fact that $B_\lambda$ takes values in $L^+_{B_{G_0}}G^{\C}_{\tau}$
and that $z\mapsto H_\lambda(z)$ is holomorphic, we deduce that
$\mu_\lambda:= (H_\lambda)^{-1}dH_\lambda$ has the desired properties. \bbox

Conversely any holomorphic potential in ${\cal H}_{-2,\infty}(L)$ produces a 
Hamiltonian stationary Lagrangian immersion as follows.
\begin{theorem}
Let $\mu_\lambda\in {\cal H}_{-2,\infty}(L)$, $p_0$ a point in $L$ and 
$H^0_\lambda$ a constant in $\Lambda G^{\C}_{\tau}$. Then
\begin{itemize}
\item there exists a unique holomorphic map $H_\lambda:L\to 
\Lambda G^{\C}_{\tau}$,
such that $dH_\lambda = H_\lambda\mu_\lambda$ and 
$H_\lambda(p_0) = H^0_\lambda$.
\item if the loop groups decomposition
$\Lambda G^{\C}_{\tau}=\Lambda G_{\tau}.L^+_{B_{G_0}}G^{\C}_{\tau}$
holds then we can apply it to $H_\lambda(z)$ for all value of $z$. It follows that there
exists two maps $F_\lambda:L\to \Lambda G_{\tau}$ and
$B_\lambda:L\to L^+_{B_{G_0}}G^{\C}_{\tau}$
such that
$$H_\lambda(z) = F_\lambda(z)B_\lambda(z),\quad \forall z\in L.$$
Then $F_\lambda$ is a lift of a (conformal) Hamiltonian stationary Lagrangian
immersion.
\end{itemize}
\end{theorem}
\emph{Proof.} Since $\mu_\lambda = \xi_\lambda dz$, with
${\partial \xi_\lambda\over \partial \overline{z}}=0$, it follows easily that
$d\mu_\lambda + \mu_\lambda\wedge \mu_\lambda = 0$, hence the existence
and the uniqueness of $H_\lambda$. Assume now that we can perform the
generalized Iwasawa decomposition $H_\lambda = F_\lambda B_\lambda$. 
It implies that
\begin{equation}\label{f-1df}
(F_\lambda)^{-1}dF_\lambda = B_\lambda\mu_\lambda(B_\lambda)^{-1}
-dB_\lambda(B_\lambda)^{-1}.
\end{equation}
Now using the fact that $\mu_\lambda\in {\cal H}_{-2,\infty}(L)$ and
$B_\lambda$ takes value in $\Lambda^+_{B_{G_0}}G^{\C}_{\tau}$,
it is easy to check that the right hand side of (\ref{f-1df}) has the form
$\sum_{k=-2}^{\infty}\hat{\alpha}_k\lambda^k$. But (\ref{f-1df}) implies also 
that this quantity should be real, i.e. a 1-form with coefficients in
$\Lambda G_{\tau}$. Hence $\alpha_\lambda:= (F_\lambda)^{-1}dF_\lambda$ reduces to
$\alpha_\lambda = \hat{\alpha}_{-2}\lambda^{-2} + \hat{\alpha}_{-1}\lambda^{-1}
+\hat{\alpha}_0 + \hat{\alpha}_1\lambda + \hat{\alpha}_2\lambda^2$ and
moreover $\hat{\alpha}_0$ is real, $\hat{\alpha}_1 = \overline{\hat{\alpha}_{-1}}$
and $\hat{\alpha}_2 = \overline{\hat{\alpha}_{-2}}$. Lastly a Taylor expansion
in $\lambda$ of (\ref{f-1df}) proves that $\hat{\alpha}_{-2}$ and $\hat{\alpha}_{-1}$
are (1,0)-forms, which ensures the result by Proposition \ref{characterization}. \bbox

        \subsection{Meromorphic potentials}

The holomorphic potentials constructed in Lemma \ref{holopotential} are 
far from being unique. Moreover they involved in general infinitely many 
holomorphic maps. These defects can be mended, provided we allow meromorphic 
potentials and under some hypotheses on $G$. We define  
$$
L_{-2,-1}\g ^{\C}_{\tau}:= \{[\lambda\mapsto \phi_\lambda]
\in \Lambda\g ^{\C}_{\tau}/
\xi_\lambda = \hat{\xi}_{-2}\lambda^{-2} + \hat{\xi}_{-1}\lambda^{-1}\}.
$$
\begin{definition}
The set of {\em meromorphic potentials}, denoted
${\cal M}_{-2,-1}(L)$, is the set of meromorphic
1-forms on $L$ with coefficients in $L_{-2,-1}\g ^{\C}_{\tau}$. So any
form $\mu_\lambda$ in ${\cal M}_{-2,-1}(L)$ has the expression
$$
\mu_\lambda = \hat{\mu}_{-2}\lambda^{-2} + \hat{\mu}_{-2}\lambda^{-2}
= (\hat{\xi}_{-2}(z)\lambda^{-2} + \hat{\xi}_{-1}(z)\lambda^{-1})dz,
$$
where $\hat{\xi}_{-2}(z)\lambda^{-2} + \hat{\xi}_{-1}(z)\lambda^{-1}\in
L_{-2,\infty}\g ^{\C}_{\tau}$.
\end{definition}
Then using the same methods as in \cite{DPW}; one can prove the following
\begin{theorem}
Assume that the conclusion of Lemma \ref{riemann-hilbert} holds. Let
$F_\lambda :L\to \Lambda G_{\tau}$ be the extended
lift of a (conformal) Hamiltonian stationary Lagrangian immersion. Then there
exists a finite subset $\{a_1,\dots ,a_p\}$ of $L$ such that
\begin{itemize}
\item there exists a holomorphic map 
$F^-_\lambda :L\setminus \{a_1,\dots ,a_p\}\to
\Lambda^-_{\star} G_{\tau}^{\C}$ and a map
$F^+_\lambda :L\setminus \{a_1,\dots ,a_p\}\to
\Lambda^+ G_{\tau}^{\C}$ such that
$$
F_\lambda(z) = F^-_\lambda(z) F^+_\lambda(z),\quad \forall
z\in L\setminus \{a_1,\dots ,a_p\}
$$
\item $z\mapsto F^-_\lambda(z)$ extends to a meromorphic map on $L$
\item the Maurer-Cartan form $\mu_{\lambda}:= (F^-_\lambda)^{-1}dF^-_\lambda$
of $F^-_\lambda$
is a meromorphic potential in ${\cal M}_{-2,-1}(L)$.
\end{itemize}
\end{theorem}
\emph{Proof.} (see \cite{DPW} for details) The decomposition
$F_\lambda(z) = F^-_\lambda(z) F^+_\lambda(z)$ is possible as soon as we can
prove that $F_\lambda(z)$ belongs to the big cell ${\cal C}$. Using Lemma
\ref{holopotential} in the same way as in \cite{DPW}, one can show that this is true
for all $z$, excepted maybe on a finite subset $\{a_1,\dots ,a_p\}\subset L$. The
second property is proved also in \cite{DPW}. The last one follows easily by
writing

$$\mu_{\lambda} = F^+_\lambda \left[ \alpha_\lambda -
(F^+_\lambda)^{-1}dF^+_\lambda \right] (F^+_\lambda)^{-1}$$
which implies on the one hand that $\mu_{\lambda}$ is in
${\cal H}_{-2,\infty}(L \setminus \{a_1,\dots ,a_p\})$, once one keep in mind the fact that
$F^+_\lambda(z)\in \Lambda^+ G_{\tau}^{\C}$.
But on the other hand
$F^-_\lambda(z)\in \Lambda^-_{\star} G_{\tau}^{\C}$
and thus there is no nonnegative power of $\lambda$ in the Fourier expansion of
$\mu_{\lambda}$. This implies the conclusion. \bbox

        \section{A list of cases}

        \subsection{The Euclidean space}

The case of \C[2] has been thoroughly studied in a first article~\cite{HR1}
including an explicit description of all tori. We will only point out
the -- obvious -- differences between \C[2] and the other Hermitian
symmetric spaces in the light of our study. Since the group of isometries
is the semi-direct product $G=U(2) \ltimes \C[2]$, we have the additional 
commutation property $[\m,\m]=0$. The equations in~(\ref{flatness})
then decouple to yield a PDE on $H$ (in $\alpha_2,\alpha_0$)
and a PDE on \m\/ with parameters $\alpha_2,\alpha_0$. Moreover the only 
nonlinearity has disappeared~! Finally, one may go even further than the 
standard analysis of that case, since $d\alpha_0 
+ \frac{1}{2}[\alpha_0 \wedge \alpha_0]=0$ implies the (local)
existence of lifts gauging $\alpha_0$ to zero (we call them spinor lifts).
At that point the problem of finding surfaces is equivalent to solving
two linear PDEs (plus the Poincar\'e integration procedure to get $f$).
But we do not need to use the coclosedness condition, and the commutation 
property is also the key point for the linear Weierstrass representation
of Lagrangian surfaces in \C[2] -- not only stationary ones.
We derive a Dirac-type equation characterizing all such surfaces 
(see~\cite{HR2}).

        \subsection{The complex projective plane} \label{CP2}

We write the projective plane as a symmetric space
\[
\C\P[2] = G/H = \frac{SU(3)}{S(U(2) \times U(1))}
\]
where\footnote{Notice that the action of $SU(3)$ is only almost 
effective, with a kernel made of three elements: the cubic roots of 
identity.}
\[
S(U(2) \times U(1)) = \left\{ \bmat{cc} A & 0 \\ 0 & \det A^{-1} \emat
, \; A \in U(2)  \right\}
\]
with Lie algebra
\[
\h = \sl(3,\C) \cap (\goth{u}(2) \oplus \goth{u}(1) )
= \left\{ \bmat{cc} X & 0 \\  0 & - \tr X \emat, \; X \in \u(2) \right\}.
\]
Here and in subsequent sections, $\tilde{X}$ will denote the 
conjugate of $X$ with respect to the real form $\g \subset \gC$; in 
the $\su(3)$ case, $\tilde{X} = - X^*$.

The quotient map is given simply by $SU(3) \to \C\P[2]$,
$g \mapsto \C g \epsilon_3$ where $\epsilon_3=(0,0,1)$.
The natural involution $\sigma$ acts on $SU(3)$ (and its differential 
on $\g = \su(3)$) by conjugation:
\[
\sigma : \bmat{cc} A & u \\ -u^* & a \emat
\mapsto \bmat{cc} \1_2 & 0 \\ 0 & -1 \emat 
\bmat{cc} A & u \\ -u^* & a \emat
\bmat{cc} \1_2 & 0 \\ 0 & -1 \emat
= \bmat{cc} A & -u \\ u^* & a \emat
\]
The Lie algebra $\g = \su(3)$ splits as the direct sum 
of the $+1$ eigenspace of $\sigma^2$, $\h$, 
and the (-1)-eigenspace 
\[
\m = \left\{ \bmat{cc} 0 & u \\ -u^* & 0 \emat , \; u \in \C[2] \right\}
\]
identified with $\C[2]$ via
\[
\bmat{cc} 0 & u \\ -u^* & 0 \emat \mapsto u  .
\]
The adjoint representation of $H$ on \m\/ is surjective and 
almost effective:
\[
\Ad \bmat{cc} A & 0 \\ 0 & \det A^{-1} \emat = \Big[ u \mapsto (\det A) A u \Big]
\]
and
\[
\ad \bmat{cc} X & 0 \\ 0 & - \tr X \emat = \Big[ u \mapsto (X+\tr X \1)u \Big]
\]
so that the complex structure is $\exp(\frac{\pi}{2}\ad Y)$ with
\[
Y = \frac{i}{3} \bmat{ccc} 1 \\ & 1 \\ && -2 \emat \; . 
\]

The order four automorphism acting on $G^{\C}$ is
\begin{equation}    \label{tauCP2}
\tau : g \mapsto  \bmat{cc|c}  & 1 &  \\ -1 & & \\ \hline & & 1 \emat
\trsp{g}^{-1} \bmat{cc|c}  & -1 &  \\ 1 & & \\ \hline & & 1 \emat .
\end{equation}
Notice that on  $G$ itself $\trsp{g}^{-1} = \bar{g}$; hence its 
differential acting on $\gC$ (still denoted $\tau$) is:
\[
\tau : X \mapsto - \bmat{cc|c}  & 1 &  \\ -1 & & \\ \hline & & 1 \emat
\trsp{X} \bmat{cc|c}  & -1 &  \\ 1 & & \\ \hline & & 1 \emat .
\]
We then have a direct sum $\gC = \gC_0 \oplus \gC_2 \oplus \gC_{-1} 
\oplus \gC_1$ with $\gC_2 = \C Y$,
\[
\g_0 = \left\{ \bmat{cc} X & 0 \\ 0 & 0 \emat ; \; X \in \su(2) \right\} 
\]
\[
\gC_{-1} = \left\{ \bmat{cc|c} && -ia \\ && b \\ \hline -ib & a \emat, 
\; a,b \in \C \right\}
\]
and
\[\gC_1  = \widetilde{\gC_{-1}}
= \left\{ \bmat{cc|c} && ia \\ && b \\ \hline ib & a \emat,\;  
a,b \in \C \right\} \; . 
\]

\begin{example}\em
\textsc{The real projective plane} \RP[2] is immersed minimally 
in \CP[2] (and its double cover is the only minimal Lagrangian sphere, 
up to unitary isometries, see~\cite{Y}). Choose the stereographic 
projection from the southern pole as  conformal coordinate chart. 
The fundamental lift (real-valued of course) is: 
\[
F(z) = \frac{1}{1+|z|^2} \bmat{ccc}
1-x^2+y^2 & -2xy & 2x \\ 
-2xy & 1+x^2-y^2 & 2y \\
-2x & -2y & 1-x^2-y^2 \emat
\]
The Maurer-Cartan form satisfies $\alpha_2 = 0$ and
\begin{eqnarray*}
(1+|z|^{2}) \alpha & = &
\underbrace{\bmat{cc|c} & -i &  \\ i & &\\ \hline &&{}\emat 
(\bar{z}dz -zd\bar{z})}_{\alpha_0}
\\
&&+\underbrace{\bmat{cc|c} && 1 \\ && -i \\ \hline -1 & i \emat 
dz}_{\alpha'_{-1}} 
+\underbrace{\bmat{cc|c} && 1 \\ && i \\ \hline -1 & -i \emat 
d\bar{z}}_{\alpha''_1}.
\end{eqnarray*}
The associated family is only a change of variable by rotation in the 
$z$-plane.
\end{example}
\begin{example}\em
\textsc{The Clifford torus} is the quotient of the standard torus 
$S^1 \times S^1 \times S^1 \subset \C[3]$ by the Hopf action; it can 
be conformally parametrized as $\C f$ where
\[
f(x+iy) = \frac{1}{\sqrt{3}}\bmat{c} e^{2ix} \\ e^{i(y\sqrt{3}-x)} \\ 
e^{-i(x+y\sqrt{3})} \emat
\]
and the fundamental lift is
\[
F(z) = \frac{1}{\sqrt{6}}\bmat{ccc}
2i e^{2ix} & 0 & \sqrt{2} e^{2ix} \\ 
-i e^{i(y\sqrt{3}-x)} & i \sqrt{3} e^{i(y\sqrt{3}-x)} & 
        \sqrt{2}e^{i(y\sqrt{3}-x)} \\ 
-i e^{-i(x+y\sqrt{3})} & -i \sqrt{3} e^{-i(x+y\sqrt{3})} 
        & \sqrt{2}e^{-i(x+y\sqrt{3})} \emat
\]
with Maurer-Cartan form
\begin{eqnarray*}
\alpha &=& 
\underbrace{ \bmat{cc|c} i & -1 \\ -1 & -i \\ \hline &&{}\emat\frac{dz}{2}
+ \bmat{cc|c} i & 1 \\ 1 & -i \\ \hline &&{}\emat\frac{d\bar{z}}{2}}_{\alpha_0}
\\
&&+ \underbrace{ \bmat{cc|c}
&& 1 \\ && -i \\ \hline -1 & i \emat \frac{dz}{\sqrt{2}}}_{\alpha'_{-1}} 
+ \underbrace{ \bmat{cc|c} && 1 \\ && i \\ \hline -1 & -i \emat 
\frac{d\bar{z}}{\sqrt{2}}}_{\alpha''_{1}} .
\end{eqnarray*}
Note that $\alpha_2 = 0$, which agrees with the fact that the Clifford torus 
is minimal.
\end{example}
\begin{example}\em
\textsc{Vacuum solutions} are obtained by taking extended lifts 
$F_\lambda = \exp(zM_\lambda+\bar{z}\widetilde{M_\lambda})$
where $M_\lambda$ is a constant in $\Lambda G_\tau$.
Equation~(\ref{flatness}) amounts to $[M_\lambda,\widetilde{M_\lambda}]=0$
For further simplification ($M_\lambda$ being constant) 
we gauge the $\gC_{-1}$ part to 
\[
\frac{e^\rho}{2}\bmat{cc|c} && 1 \\ && -i \\ \hline -1 & i \emat .
\]
That yields a family parametrized by complex numbers $b,c$ such that
$e^{2\rho} = 8 \Im (\bar{b}c) >0$
\[
M_\lambda = \bmat{ccc}
b - i \lambda^{-2} a  & c & \lambda^{-1} \frac{e^\rho}{2} 
\\
c   & -b - i \lambda^{-2} a  & -i \lambda^{-1} \frac{e^\rho}{2} 
\\ 
-\lambda^{-1} \frac{e^\rho}{2} & i\lambda^{-1} \frac{e^\rho}{2} & 
2\lambda^{-2}ia
\emat , \quad a = -\frac{\bar{c}+i\bar{b}}{3}
\]
Minimal conformal immersions correspond to $a=0$ so
$c=ib$ and $e^{\rho} = 2\sqrt{2}|b|$
\[
M_\lambda = \bmat{ccc}
b & ib & \lambda^{-1} \sqrt{2} |b| 
\\
ib   & -b  & -i \lambda^{-1} \sqrt{2} |b| 
\\ 
-\lambda^{-1} \sqrt{2} |b| & i\lambda^{-1} \sqrt{2} |b| & 0
\emat
\]
and for $b = \frac{i}{2}$ we recognize the Clifford torus above.
Variations in $|b|$ amount to trivial scale changes, but changes in 
the argument of $b$ yield different examples (not gauge-equivalent);
these however may not be periodic.
\end{example}

        \subsection{Lagrangian cones in \C[3]}  \label{cones_frames}
        
We explain here how a similar formalism applies as well to Lagrangian 
cones in \C[3], and how this relates formally to the previous section 
(knowing that such cones are intimately associated to Lagrangian surfaces 
in \CP[2] as explained in \S\ref{cones}). This association is not 
one to one, but to each Lagrangian surface corresponds exactly
a circle of cones, namely the orbit under the Hopf action of any member.
To make this relation visually explicit, we will overline with a $\check{}$ 
the corresponding quantity in $S^5$; for instance, the Legendrian map 
$\check{f} : L \to S^5$ projects down to a Lagrangian map $f : L \to 
\CP[2]$.

We view now $S^5$ as the reductive\footnote{but not symmetric.} space 
$\check{G}/\check{H} = U(3)/\check{H}$ where
\[
\check{H} = \left\{ \bmat{cc} A & 0 \\ 0 & 1 \emat,\; A \in U(2) \right\}
\]
is the isotropy group of $\epsilon_3$. The quotient map is 
$g \mapsto g \epsilon_3$. We have the reductive splitting
$\check{\g} = \u(3) = \check{\h} \oplus \check{\m}$ with
\[
\check{\h} = \left\{ \bmat{cc} X & 0 \\ 0 & 0 \emat,\; X \in \u(2) \right\}
\; , \;
\check{\m} = \left\{ \bmat{cc} 0 & u \\ -u^*  & ia \emat , \; 
a\in\R, \; u \in \C[2] \right\}.
\]
The same order four automorphism $\tau$ (formally) as in the previous section 
(see formula~(\ref{tauCP2})) acts on $U(3)$ and splits the Lie algebra 
$\check{\g}^{\C} = \goth{gl}(3,\C)$ into four eigenspaces: 
$\check{\g}^{\C}_0 \oplus \check{\g}^{\C}_2 \oplus \check{\g}^{\C}_{-1} 
\oplus \check{\g}^{\C}_1$:
\[ 
\check{\g}_2 = \R \check{Y} \oplus \R \check{Z}  \textrm{ with }
\check{Y} = \bmat{ccc} i \\ & i \\ & & 0 \emat , \;
\check{Z} = \bmat{ccc} 0 \\ & 0 \\ & & i \emat
\]
and $\check{\g}_0 = \g_0$, $\check{\g}^{\C}_{-1} = \gC_{-1}$, 
$\check{\g}^{\C}_{1} = \gC_{1}$. Comparing with \S\ref{CP2}, the 
only differences are (i) $\check{\g}_2$ is two dimensional and 
(ii) the complex structure changes from $Y$ to $\check{Y}$. Nota bene: 
the contact distribution is generated as the orbit under $U(3)$ of 
the subspace 
\[
\check{\p} = \left\{ \bmat{cc} 0 & u \\ -u^*  & 0 \emat , \; 
u \in \C[2] \right\} \simeq \C[2],
\]
endowed with the complex structure $\ad \check{Y}$. Obviously 
$\check{\p}^{\C} = \check{\g}^{\C}_{-1} \oplus \check{\g}^{\C}_1$
and $\check{\m} = \check{\p} \oplus \R \check{Z}$.

Consider a Lagrangian cone $\mathcal{C}$; its link is Legendrian, namely 
satisfies that: (i) its tangent bundle lies in the 
contact distribution $\Pi$ ($T_xM \perp ix$) and (ii) the tangent space
$T_x M$ is Lagrangian in $\Pi_x$. 
Letting $\check{f} : L \to S^5$ be a conformal parametrization of $M$, 
conditions (i) and (ii) above amount to the existence of a (unique) 
fundamental lift $\check{F} \in U(3)$ such that:
\begin{equation}        \label{conelagconf}
d\check{f}_z = e^{\rho(z)} \check{F}(z)(\epsilon_1 dx + \epsilon_2 dy)
\end{equation}
which can be rewritten in terms of $\check{\alpha} = \check{F}^{-1} d\check{F}$
\begin{equation}
\check{\alpha} \epsilon_3 = e^{\rho} (\epsilon_1 dx + \epsilon_2 dy)
= e^{\rho} (\epsilon dz + \tilde{\epsilon} d\bar{z})
\end{equation}
where
\[
\epsilon_1 = \bmat{cc|c} && 1 \\ && 0 \\ \hline -1 & 0 \emat , \;
\epsilon_2 = \bmat{cc|c} && 0 \\ && 1 \\ \hline 0 & -1 \emat,
\]\[
\epsilon = \frac{\epsilon_1-i\epsilon_2}{2} 
= \frac{1}{2} \bmat{cc|c} && 1 \\ && -i \\ \hline -1 & i \emat, \;
\tilde{\epsilon} = \frac{1}{2} \bmat{cc|c} && 1\\ && i\\ \hline -1 & -i \emat.
\]
Recall that $\det_{\C} \check{F} = e^{i\beta}$ where $\beta$ is the 
Lagrangian angle. Extending to the orbit under the gauge action of
$C^\infty(L,\check{G}_0)$ with $\check{G}_0 = \check{H} \cap SU(3) 
\simeq SU(2)$, we define Lagrangian lifts by the property that
$\check{\alpha}_{\m} = \check{\alpha}'_{-1} + \check{\alpha}''_1$.
Notice that the condition is more complicated here because we 
need to assume that the map is horizontal (i.e. lies in the contact 
distribution), which excludes components along $\check{Z}$. 

We can now characterize Hamiltonian stationary Lagrangian cones 
either intrinsically through the following
\begin{theorem}
A conical Lagrangian singularity (whose intersection with $S^5$ is
conformally parametrized) is exactly obtained by integrating a flat 
$\u(3)$-valued 1-form $\omega$ satisfying $\check{\alpha}_\m = 
\check{\alpha}'_{-1}+\check{\alpha}''_1$ (hence $\check{\alpha}_2$ 
lies in $\C \check{Y}$ since $\R \check{Z} = \check{\m} \cap \check{\g}_2$).
Furthermore, the immersion is H-minimal (resp. special Lagrangian) if
$\alpha_2$ is coclosed (resp. vanishes).
\end{theorem}
Denoting $\Lambda \u(3)_{\tau,\p}$ the subspace of the twisted loop-algebra 
where $Z$-part vanishes, there is an interesting bijective correspondence 
between this subspace and the loop algebra $\Lambda \su(3)_\tau$, mapping 
flat extended connection forms to flat extended connections forms, which 
leaves all matrix coefficients unchanged but for the $\check{\g}_2$ part where 
the complex structure $\check{Y}$ is mapped to complex structure~$Y$:
\[ 
\bmat{ccc} a \\ & b \\ && 0 \emat \mapsto 
\frac{1}{3}\bmat{ccc} 2a-b \\ & 2b-a \\ 
&&0 \emat .
\]

Or one can associate to the cone the projected Lagrangian surface in 
$\CP^2$ with the following data: a map $f=\C \check{f} : L \to \CP[2]$ 
with a Lagrangian lift $F$. We claim that $F = e^{i\beta/3} \check{F}$ 
is such a lift. Indeed $F$ lifts $f$ since $\C F \epsilon_3 = f$ and 
$\det F=1$. To prove that $F$ is Lagrangian consider its Maurer-Cartan form
\[
\alpha = \check{\alpha} + \frac{d\beta}{3} \1
= \check{\alpha}_0 + \check{\alpha}'_{-1} + \check{\alpha}''_1
+ \underbrace{\frac{d\beta}{2} \check{Y} - \frac{id\beta}{3}\1}_{
 \displaystyle \frac{d\beta}{2}Y} .
\]
Obviously $\alpha''_{-1} = 0$. Furthermore we see that the Lagrangian 
angle of the cone is equal to the Lagrangian angle of the surface in 
$\CP[2]$. It may be noted that the fundamental lift is mapped thus to 
the fundamental lift.

        \subsection{The complex hyperbolic plane}

The non compact dual of $\CP[2]$ is the complex hyperbolic space
\[
\C\D[2]=SU(2,1)/S(U(2)\times U(1))
\]
where
\[
SU(2,1) = \left\{ g \in SL(3,\C), \;
g B g^* = B = \bmat{ccc} 1 \\ & 1 \\ & & -1 \emat
\right\} 
\]
with Lie algebra
\[
\su(2,1) = \left\{ \bmat{cc} X & v \\ v^* & -\tr X \emat , \; 
X \in \u(2), \; v \in \C[2] \right\}.
\]
The same automorphism $\tau$ acts on $M$.

        \subsection{$\CP[1] \times \CP[1]$}

We consider the following Hermitian symmetric space:
\[
\CP[1] \times \CP[1] = \frac{SU(2) \times SU(2)}{U(1) \times U(1)} 
\]
($G = SU(2) \times SU(2)$ will be written as bloc-diagonal four by 
four matrices) and define an order four automorphism 
$\tau : g \mapsto TgT^{-1}$ where 
\[ 
T = \bmat{cc|cc} 
 &&& -1 \\ && 1 \\ \hline & 1 \\ 1 &  \emat .
\] 
Its differential at identity diagonalizes on $\gC = \sl(2) \oplus \sl(2)$ 
with eigenspaces $\gC_0 = \C X$, $\gC_2 = \C Y$ (as usual $Y$ is the 
complex structure)
\[
X = \frac{1}{2} \bmat{cc|cc}
-i&  \\  & i \\ \hline && i  \\ &&& -i \emat ,
\qquad
Y =\frac{1}{2} \bmat{cc|cc}
-i&  \\  & i \\ \hline && -i  \\ &&& i \emat ,
\]
\[
\gC_{-1} = \left\{ \bmat{cc|cc}
 & v \\ u & \\ \hline &&& iu \\ && iv \emat u,v \in \C \right\} ,
\]\[
\gC_{1} = \left\{ \bmat{cc|cc}
 & v \\ u&  \\ \hline &&& -iu \\ && -iv \emat u,v \in \C \right\} .
\]

        \subsection{The non compact dual of $\CP[1] \times \CP[1]$}

As expected, the situation is very close to its compact dual.
\[ 
M = G/H = \frac{SU(1,1) \times SU(1,1)}{U(1) \times U(1)}
\]
the automorphism  has the same expression and so do the eigenspaces.

\lf

\lf\lf\lf
{ 
Fr\'ed\'eric H\'elein \lf
\texttt{helein@cmla.ens-cachan.fr} \lf
\lf
Pascal Romon \lf
\texttt{romon@cmla.ens-cachan.fr}
}

\begin{thebibliography}{XXXXX}
    
\bibitem[A]{A} Henri Anciaux:
\emph{An isoperimetric inequality for Hamiltonian stationary Lagrangian 
tori in $\C[2]$ related to Oh's conjecture}, CMLA preprint.

\bibitem[Br]{Br} Robert Bryant:
\emph{Minimal Lagrangian submanifolds of K\"ahler-Einstein manifolds},
Differential geometry and differential equations (Shanghai, 1985), 1-12, 
Lecture Notes in Math., 1255, 
Springer, Berlin-New York, 1987.

\bibitem[BP]{BP} Francis E. Burstall, Franz Pedit:
\emph{Dressing orbits of harmonic maps}, 
Duke Math. J. 80 (1995), nr 2, 353-382.

\bibitem[BRa]{BRa} Francis E. Burstall, John H. Rawnsley:
\emph{Twistor theory for Riemannian symmetric spaces},
Lecture Notes in Mathematics 1424, Springer-Verlag (1990).

\bibitem[CaUr]{CaUr} Ildefonso Castro, Francisco Urbano, 
\emph{Examples of unstable Hamiltonian-minimal Lagrangian tori in $\C^2$}, 
Composition Mathematica 111 (1998), 1-14.

\bibitem[CM]{CM} Bang-Yen Chen, Jean-Marie Morvan: 
{\em Deformations of isotropic submanifolds in K\"ahler manifolds}, 
J. of Geometry and Physics 13 (1994), 79-104.

\bibitem[DPW]{DPW} Josef Dorfmeister, Franz Pedit, Hong-You Wu: 
{\em Weierstrass type representation of harmonic maps into symmetric spaces}, 
Comm. in Analysis and Geometry Vol 6, Number 4 (1998), 633-668.

\bibitem[GaHuLa]{GaHuLa} 
Sylvestre Gallot, Dominique Hulin, Jacques Lafontaine: 
\emph{Riemannian geometry}, 
Universitext, Springer-Verlag. XI,(1987).

\bibitem[G]{G} Phillip Griffiths: 
\emph{On Cartan's method of Lie groups and moving frames as applied 
to uniqueness and existence questions in differential geometry},
Duke Math. J. 41, 775-814 (1974).

\bibitem[HaL]{HaL} Reese Harvey, H. Blaine Lawson: 
{\em Calibrated geometries}, 
Acta Mathematica 148 (1982), 47-157.

\bibitem[Has]{Has} Mark Haskins:
\emph{Special Lagrangian Cones}, \textsf{arXiv:math.DG/0005164}.

\bibitem[H]{H} Fr\'ed\'eric H\'elein:
\emph{Willmore surfaces and loop groups}, 
J. Differential Geometry 50 (1998), p. 331-385.

\bibitem[HR1]{HR1} Fr\'ed\'eric H\'elein, Pascal Romon:
\emph{Hamiltonian stationary Lagrangian surfaces in \C[2]},
\textsf{arXiv:math.DG/0009202}.

\bibitem[HR2]{HR2} Fr\'ed\'eric H\'elein, Pascal Romon:
\emph{Weierstrass representation of Lagrangian surfaces in 
four-dimensional space using spinors and quaternions},
to appear in Commentarii Mathematici Helvetici.

\bibitem[Hi]{Hi} Masanori Higaki: 
\emph{Actions of loop groups on the space of harmonic maps into 
reductive homogeneous spaces},
J. Math. Sci. Univ. Tokyo 5 (1998), 401-421.

\bibitem[O]{O} Yong-Geun Oh, 
\emph{Volume minimization of Lagrangian submanifolds 
under Hamiltonian deformations}, Math. Z. 212 (1993), 175-192.

\bibitem[PrS]{PrS} Andrew Pressley, Graeme Segal:
{\em Loop groups}, 
Oxford Mathematical Monographs, Clarendon Press, Oxford 1986.

\bibitem[Re1]{Re1} Helmut Reckziegel:
\emph{Horizontal lifts of isometric immersions into the bundle space 
of a pseudo-Riemannian submersion}, 
Global differential geometry and global analysis (1984), 
Lecture Notes in Mathematics 1156, Springer Verlag (1985), 264-279.

\bibitem[Re2]{Re2} Helmut Reckziegel:
\emph{A correspondence between horizontal submanifolds of Sasakian 
manifolds and totally real submanifolds of K\"ahlerian submanifolds},
Pap. Colloq., Hajduszoboszlo/Hung. 1984, Vol. 2, Colloq. Math. Soc. 
Janos Bolyai 46, 1063-1081 (1988).

\bibitem[SaSa] {SaSa} Mikio Sato, Yasuko Sato:
{\em Soliton equations as dynamical systems in infinite dimensional
Grassmann manifolds}, 
in Nonlinear partial differential equations in applied science 
(Tokyo, 1982), 259--271, North-Holland Math. Stud., 81, (1983).

\bibitem[ScWo]{ScWo} Richard Schoen, Jon Wolfson,
\emph{Minimizing volume among Lagrangian submanifolds}, 
Proc. Sympos. Pure Math., 65, Amer. Math. Soc. (1999).

\bibitem[SW] {SW} Graeme Segal, George Wilson:
{\em Loop groups and equations of KdV type}, Publ. Math. IHES 61(1985), 5-65,
also published in {\em Surveys in differential geometry: integrable systems},
Vol. IV, C.L. Terng, K. Uhlenbeck editors, International Press 1998.

\bibitem[U]{U} Karen Uhlenbeck:
{\em Harmonic maps into Lie groups}, 
Journal of Differential Geometry 30 (1989), 1-50.

\bibitem[Wo]{Wo} Jon Wolfson: 
{\em Minimal Lagrangian diffeomorphisms and the Monge-Amp\`ere equations}, 
J. Differential Geometry 46 (1997), 335-373.

\bibitem[Y]{Y} Shing-Tung Yau:
\emph{Submanifolds with constant mean curvature I,II},
Amer. J. Math., 96 (1974), 346-366; ibid. 96 (1975), 76-100.
\end{thebibliography}
\end{document}